\documentclass[a4paper,leqno,11pt]{amsart}
\usepackage{amsfonts,amssymb,verbatim,amsmath,amsthm,latexsym,textcomp,amscd}
\usepackage{latexsym,amsfonts,amssymb,epsfig,verbatim}
\usepackage{amsmath,amsthm,amssymb,latexsym,graphics,textcomp}
\usepackage{graphicx}
\usepackage{color}
\usepackage{url}
\usepackage{enumerate}
\usepackage[mathscr]{euscript}
\usepackage{tikz}
\usepackage{dsfont}
\usetikzlibrary{matrix}
\usepackage{sseq}
\usepackage{hyperref}

\input xy
\xyoption{all}

\theoremstyle{plain}

\theoremstyle{definition}
\newtheorem{theorem}{Theorem}[section]
\newtheorem{prop}[theorem]{Proposition}
\newtheorem{lemma}[theorem]{Lemma}

\newtheorem{subsec}[theorem]{}

\newtheorem{thm}[theorem]{Theorem}

\theoremstyle{remark}
\swapnumbers

\newcommand{\C}{{\mathbb C}}

\newcommand{\Z}{{\mathbb{Z}}}
\newcommand{\R}{{\mathbb R}}

\newcommand\FF{{\mathcal F}}

\newcommand\LL{{\mathcal L}}
\newcommand\MM{{\mathcal M}}

\newcommand\PP{{\mathcal P}}

\newcommand\PMF{{\PP\kern-2pt\MM\FF}}

\newcommand\PML{{\PP\kern-2pt\MM\LL}}

\newcommand{\fsubd}{\mathrel{{\scriptstyle\searrow}\kern-1ex^d\kern0.5ex}}
\newcommand{\bsubd}{\mathrel{{\scriptstyle\swarrow}\kern-1.6ex^d\kern0.8ex}}
\newcommand{\fsubeq}{\mathrel{\raise-.7ex\hbox{$\overset{\searrow}{=}$}}}
\newcommand{\bsubeq}{\mathrel{\raise-.7ex\hbox{$\overset{\swarrow}{=}$}}}

\newcommand{\tsh}[1]{\left\{\kern-.9ex\left\{#1\right\}\kern-.9ex\right\}}

\makeatletter
\@namedef{subjclassname@2020}{%
  \textup{2020} Mathematics Subject Classification}
\makeatother

\title{Characteristic classes on certain quotients of Stiefel manifolds}
\author{Debanil Dasgupta}

\email{debanil12@gmail.com;}
\address{Stat-Math Unit,
	Indian Statistical Institute,
	B. T. Road, Kolkata-700108, India.}

\subjclass[2020]{Primary: 57R20, 57N65} 
\keywords{Stiefel manifolds, characteristic classes, embedding.}

\begin{document}

\maketitle

\begin{abstract}
	In this paper, we calculate characteristic classes for certain quotients of real Stiefel manifolds $V_{n,k}$ and then derive results on certain numerical invariants, such as characteristic rank and skew embedding dimension, for those spaces .  
\end{abstract}

\section{Introduction}

This paper explores certain topological results on quotients of Stiefel manifolds. The perspective behind these results are cohomology calculations which in turn lead to geometric consequences. The quotients of Stiefel manifolds form a nice collection of homogeneous spaces which are amenable to computational techniques. We consider
\[ V_{n,k} = \mbox{ real Stiefel manifold of } k\mbox{-orthonormal vectors in } \R^n \cong O(n)/O(n-k), \]
\[ W_{n,k} = \mbox{ complex Stiefel manifold of } k\mbox{-orthonormal vectors in } \C^n \cong U(n)/U(n-k). \]
These manifolds are ubiquitous in the literature, and many of their features have been extensively studied. 

The Stiefel manifolds have an action by the orthogonal group $O(k)$ in the real case, and $U(k)$ in the complex case, whose quotients are the Grassmannian manifolds. The Grassmannians are in a sense the quotient of a Lie group by a parabolic subgroup. One may instead take the quotient by a Borel subgroup and end up with flag manifolds. These manifolds are also widely studied in topology and algebraic geometry. 

A different viewpoint is provided when we consider actions by cyclic subgroups of the orthogonal group. As an example, one has the action of $S^1$ on $W_{n,k}$ by diagonal matrices in $U(k)$. The quotient space is called the projective Stiefel manifold $PW_{n,k}$. The real version of this is the quotient $PV_{n,k}= V_{n,k}/C_2$, where the cyclic group $C_2$ of order $2$ acts via $\pm 1$. One may also consider a cyclic group action on $W_{n,k}$ by $m^{th}$ roots of unity, and call the quotient $W_{n,k;m}$. 

The Stiefel manifolds are parallelizable except for the spheres. For the quotients of Stiefel manifolds one has a nice approach towards calculation of their tangent bundle \cite{Lam1975}. There is also a direct method to compute their cohomology via the Serre spectral sequence \cite{Basu2021}. This allows us to explicitly calculate characteristic classes for these manifolds, such as the Stiefel Whitney classes, or the Pontrjagin classes. These computations may then be used to deduce non-embedding and non-immersion results into Euclidean spaces. Further, we may explore the question whether these manifolds are parallelizable, and obtain bounds on the number of linearly independent vector fields. Partial answers to these questions are provided in \cite{Sankaran2006}. 

One may also consider skew embeddings of these manifolds, that is, embeddings in which the affine spaces corresponding to the tangent spaces in the embedding are skew. This is related to an embedding of tangent bundle of the ordered configuration space. Bounds on such embeddings are also related to Stiefel Whitney classes \cite{Baralic2010}. This allows us to compute bounds on skew embeddings of $PV_{n,k}$ \ref{skewPV}. 

Recently in \cite{Korbas2016}, in an attempt to compute the cohomology of the oriented Grassmannian, the authors considered the question of representing cohomology classes via Stiefel Whitney classes of bundles. More precisely, we consider the subalgebra of $H^\ast(X;\Z_2)$ generated by Stiefel Whitney classes of vector bundles. The question posed for a space $X$ is the highest degree $d$ such that every class of degree $\leq d$ lies in this subalgebra. Such a $d$ is called the ucharrank$(X)$. This was computed for the Stiefel manifolds in \cite{Korbas2012}. We compute the same for $PV_{n,k}$ in a significant number of  cases \ref{charrkPV}. 

In this paper, we also consider the diagonal inclusion of $S^1$ inside a maximal torus of $O(k)$, and quotient the Stiefel manifold with the associated action. This is the space $Y_{n,k} = V_{n,2k}/S^1$ with the action via $S^1\cong SO(2)$ included in $O(2k)$ as block diagonal. The manifolds $Y_{n,k}$ are called circle quotient Stiefel manifolds. We compute it's tangent bundle \ref{tgtbdleY},  and show that it is not parallelizable if $n-2k\geq 4$, \eqref{Yparallel}. Further, we compute it's cohomology \ref{cohY}, skew embedding dimension, and upper characteristic rank .  

 \subsection{Organization}
In section $2$ we deal with skew embedding and characteristic rank of projective Stiefel manifolds. In section $3$ we define circle quotient Stiefel manifolds, obtain a description of its tangent bundle and calculate its $\Z_2$-cohomology. Finally in chapter $4$ we discuss the parallelizability, skew embedding, characteristic rank of circle quotient Stiefel manifolds.

\section{Computations for projective Stiefel manifolds}
In this section, we describe some geometric consequences for the real projective Stiefel manifolds, which are derivable from the cohomology. We only use the cohomology with $\Z_2$-coefficients in this section. Computations with $\Z$-coefficients or $\Z_{(2)}$-coefficients is more involved even in the case of Stiefel manifolds particularly when $n< 2k$. 

Recall that the $\Z_2$-cohomology of $PV_{n,k}$ was determined in \cite{Gitler_Handel1968}.  
$$ H^*(PV_{n,k};\Z_2)= \Z_2[x]/(x^N) \otimes V(A) , \text{ for } k<n ,$$ 
where  $N= \text{min}\{ j \mid n-k<j \leq n, \text{ and } \binom{n}{j} \text{ is odd } \} $, $ A= \{ y_j \mid n-k \leq j <n \} - \{ y_{_{N-1}}\} $ and $ |x|=1, |y_j|=j $. This is computed using the Serre spectral sequence for the fibration $V_{n,k} \to PV_{n,k}\to \R P^\infty$. The differentials are computed using the commutative diagram of fibrations 
\[ \xymatrix{ V_{n,k} \ar[r] \ar[d] & V_{n,k} \ar[d] \\ 
                   PV_{n,k} \ar[r]\ar[d] & BO(n-k) \ar[d]\\ 
                   \R P^\infty \ar[r] & BO(n).} \]
In fact, the bottom square is a homotopy pullback diagram, which may be used to describe homotopy classes of maps into $PV_{n,k}$. Proceeding in this direction, one discovers that the $PV_{n,k}$ classifies line bundles $L$ such that $nL$ has $k$ linearly independent sections.  

The tangent bundle of $PV_{n,k} $ was determined in \cite{Lam1975} and satisfies the following relation 
\begin{equation}\label{tgtpvnk}
 T(PV_{n,k}) \oplus \binom{k+1}{2}\epsilon = nk \zeta_{n,k}, 
\end{equation}
where $\zeta_{n,k}$ is the Hopf bundle over $PV_{n,k}$.  The computation of the tangent bundle is done via the $2^{k-1}$-sheeted covering space $PV_{n,k} \to F_{k}(\R^n)$, where the latter is the space of flags $V_0\subset V_1\subset \cdots \subset V_k$ of subspaces of $\R^n$ with $\dim(V_i)=i$. 

The computation of the tangent bundle \eqref{tgtpvnk} allows us to calculate the Stiefel Whitney classes of $PV_{n,k}$. These Stiefel-Whitney classes may also be calculated using the cohomology and Steenrod operations via Wu's formula. Note that the total Stiefel Whitney class of $\zeta_{n,k}$ is described by 
\[ w(\zeta_{n,k})=1+x.\]
By the Whitney sum formula, we obtain, 
\begin{equation} \label{stfwhitpvnk}
 w(T(PV_{n,k}))= (1+x)^{nk}.
\end{equation}

\subsection{Skew embeddings of $PV_{n,k}$}
The  Stiefel Whitney classes for a manifold may also be viewed as obstructions to trivializing vector bundles, and also to constructing linearly independent sections therein. The span of a manifold is the maximum number of linearly independent sections. For the projective Stiefel manifold, there are bounds on the span proved using the Stiefel Whitney classes and also by $K$-theory calculations \cite{Sankaran2006}. 

A general embedding problem for a manifold $M^n$ of dimension $n$ seeks to find the precise $k$ such that $M^n$ embeds in $\R^{n+k}$. An analogous statement may also be formulated for immersions. In this case, the Stiefel Whitney class of the stable normal bundle provides an obstruction to the immersion dimension. More precisely, consider $\bar{w}(M)  = w(M)^{-1}$, and suppose $\bar{w}_k(M)\neq 0$. Then, the manifold does not immerse in $\R^{n+k-1}$.   

An embedding of a manifold inside the Euclidean space $\R^N$ is called totally skew if the affine subspaces of $\R^N$ associated to the tangent space at different points are skew.  For a smooth manifold $M$, we define 
$$N(M)= \text{min} \{ n \mid M \text{ admits a skew embedding in }\R^n \} $$
Suppose  an $n$ dimensional manifold $M$ admits a skew embedding inside $R^N$. Then the ``skewness" condition ensures that there is a vector bundle monomorphism $$ T(F_2(M))\oplus \epsilon \longrightarrow F_2(M)\times R^N, $$ where $F_2(M)=M\times M - \Delta(M) $, and this in turns implies that if $\bar{w}_k(T(F_2(M)) \neq 0$ then $N \geq 2N+k+1$. In particular, this gives us a lower bound for $N(M)$ : $$ N(M)\geq 2n+k+1 .$$
The authors in \cite{Baralic2010} then found a condition in terms of the Stiefel Whitney classes of $M$ to produce a $k$ for which $\bar{w}_k(T(F_2(M)) \neq 0$.
\begin{thm}\cite{Baralic2010}\label{skewCond}
If $k:= \text{max}\{ i \mid \bar{w}_i(M) \neq 0 \}$, then $\bar{w}_{2k}(T(F_2(M)) \neq 0$ and hence $ N(M)\geq 2n+2k+1 .$
\end{thm}
Then we have a direct consequence of the above result for $M=PV_{n,k}$.
\begin{thm}\label{skewPV}
$N(PV_{n,k}) $ satisfies the following inequality 
$$N(PV_{n,k}) \geq 2\dim(PV_{n,k})+2m+1,$$
where $ m= \text{max} \{j \mid\binom{nk+j-1}{nk-1}\neq 0  \text{ and } 0\leqslant j \leqslant N-1 \}$.
\end{thm}
\begin{proof}
 
From the description of  the stable tangent bundle of $PV _{n,k}$, we obtain the total Stiefel-Whitney class of $PV_{n,k}$ as follows $$ w(PV_{n,k})= (1+x)^{nk}.$$ From this we can see $$ m= \text{max} \{j \mid\binom{nk+j-1}{nk-1}\neq 0  \text{ and } 0\leqslant j \leqslant N-1 \}= \text{max} \{ j\mid \bar{w}_j(PV_{n,k}) \neq 0\}.$$\par
Applying \ref{skewCond} we get $N(PV_{n,k}) \geqslant 2\dim(PV_{n,k})+2m+1$.
\end{proof}

\subsection{Characteristic rank of projective Stiefel manifolds}
Let $X$ be a connected finite CW-complex and $\xi$ a real vector bundle over X.  The \textbf{characteristic rank} of $\xi$ over $X$, denoted by $charrank_X(\xi)$, is by definition the largest integer $k, 0\leqslant k\leqslant dim(X)$, such that every cohomology class $x\in H^j(X;\Z/2),0\leqslant j\leqslant k$, is a polynomial in the Stiefel-Whitney classes $w_i(\xi)$.The \textbf{upper characteristic rank} of $X$, denoted by $ucharrank(X)$, is the maximum of $charrank_X(\xi)$ as $\xi$ varies over all vector bundles over $X$.\par
The main result regarding the characteristic rank for Stiefel manifolds shows that the lowest degree non-zero class in cohomology usually does not arise as a Stiefel Whitney class of a vector bundle, \cite{Korbas2012}. On the other hand, when we consider a quotient of a Stiefel manifold $V_{n,k}/G$, the classes in low degrees naturally arise from characteristic classes of $G$-representations. In the case of the projective Stiefel manifold with $\Z/2$-coefficients, the class $x$ and it's powers are expressable in terms of Stiefel Whitney classes of bundles. The remaining part which is additively an exterior algebra pulls back non-trivially to the Stiefel manifold, and this is usually not representable in terms of Stiefel Whitney classes of bundles. We elaborate this in the next few results. 

We start with a lemma which will be used later.  
\begin{lemma}\label{lempvnk}
For any real vector bundle $\xi$ over $V_{8r+1,2}$ with $r>1$, the class  $w_{8r}(\xi)=0$.
\end{lemma}
\begin{proof}
We have the following pushout diagram $$\xymatrix{\R P^{n-2}\times V_{n-1,k-1} \ar[r]\ar[d] & V_{n-1,k-1}\ar[d] \\ \R P^{n-1}\times V_{n-1,k-1}\ar[r] & V_{n,k}  }$$ This tells us that the cofiber of the map $V_{n-1,k-1} \longrightarrow V_{n,k}$ is $\Sigma^{n-1}(V_{n-1,k-1})_+$. \par
Now taking $(n,k)= (8r+1,2)$ gives the following cofiber sequence $$ S^{8r-1}\longrightarrow V_{8r+1,2} \longrightarrow \Sigma^{8r}S^{8r-1}_+ .$$Applying $\widetilde{KO} $ on it we get the following exact sequence $$ \widetilde{KO}(\Sigma^{8r}S^{8r-1}_+)\longrightarrow \widetilde{KO}(V_{8r+1,2}) \longrightarrow \widetilde{KO}(S^{8r-1}).$$ Now we know $\widetilde{KO}(S^{8r-1}) = \pi_{8r-1}(BO) = \pi_7(BO) = 0$. So $\widetilde{KO}(\Sigma^{8r}S^{8r-1}_+) \longrightarrow \widetilde{KO}(V_{8r+1,2})$ must be a surjection.  But we have 
\begin{align*}
  \widetilde{KO}(\Sigma^{8r}S^{8r-1}_+) &= \widetilde{KO}(S^{16r-1}\vee S^{8r} )\\
    &= \widetilde{KO}(S^{16r-1})\oplus \widetilde{KO}(S^{8r} )\\ 
    &= \pi_7(BO)\oplus \widetilde{KO}(S^{8r} )\\ 
    &= \widetilde{KO}(S^{8r} )
    \end{align*}
So we have a surjection $\widetilde{KO}(S^{8r}) \longrightarrow \widetilde{KO}(V_{8r+1,2}) $. This means any stable bundle over $V_{8r+1,2}$ is a pullback of a stable bundle over $S^{8r}$ by this composed map. But $w_n(\xi)=0$
for any vector bundle $\xi$ over $S^n$ whenever $n \neq 1,2,4,8$ \cite{Milnor1958}. So for $r>1$, $w_{8r}(\xi)=0$ for any vector bundle $\xi$ over $S^{8r}$ and hence same is true for $V_{8r+1,2} $ as well.
\end{proof}

\begin{thm}\label{charrkPV}
If $n-k=5$, $6$ or $\geq 9$, the upper characteristic rank of $PV_{n,k}$ is given by 
$$\mbox{ucharrank}(PV_{n,k})=\begin{cases} n-k-1 & \mbox{if } 2\mid \binom{n}{k-1} \\ 
n-k &\mbox{if } 2\nmid \binom{n}{k-1}.\end{cases}$$ 
\end{thm}
\begin{proof}
First we consider the case $n-k \neq N-1$, then $\text{ucharrank}(PV_{n,k}) < n-k$, since $y_{n-k} \in H^{n-k}(PV_{n,k})$ can not be a polynomial combination of Stiefel-Whitney classes of some bundle $\xi$ over $PV_{n,k}$. Because otherwise the pullback of $y_{n-k}$ via the map $V_{n,k} \longrightarrow PV_{n,k}$ would be a polynomial combination of Stiefel-Whitney classes of the pullback bundle of $\xi$ over $V_{n,k}$ contradicting \cite{Korbas2012}. Moreover we can easily see $\text{charrank}(\zeta_{n,k})= n-k-1$. So in this case $\text{ucharrank}(PV_{n,k})=n-k-1$.  \par
 On the other hand for $n-k = N-1$, $H^{n-k}(PV_{n,k})= \Z_2\langle x^{n-k}\rangle$ and  $H^{n-k+1}(PV_{n,k})= \Z_2\langle  y_{n-k+1}\rangle$. If $n-k+1$ is not a power of $2$, then $y_{n-k+1}$ cannot be a polynomial combination of Stiefel-Whitney classes because of the cohomology ring structure of $H^*(PV_{n,k};\Z_2)$ and from the fact that the Stiefel-Whitney classes of a bundle are generated by the Stiefel-Whitney classes of degree powers of $2$ of that bundle over the Steenrod algebra (this fact is a consequence of Wu's formula) \cite{MilnorStasheff1974}. So in this case $\text{ucharrank}(PV_{n,k})=n-k$ since $\text{charrank}(\zeta_{n,k})=n-k$. \par
 Now we consider the case where $n-k = N-1$, and $n-k+1=2^s$, for some $s>3$.  In this case the only way $y_{n-k+1}$ can be a polynomial in Stiefel-Whitney classes is if $y_{n-k+1}= w_{n-k+1} $. But the lemma \ref{lempvnk} rules out this possibility and we get $\text{ucharrank}(PV_{n,k})=n-k$.
\end{proof}

We also note that one may carry the argument above forward to even obtain partial results in the cases $n-k=1,2,3,4,7,8$. For example, if $n-k=1$, we observe that if $n\equiv 2,3 \pmod{4}$, then $N=2$, so the first class in the exterior algebra part is $y_2$. This pulls back non-trivially to $V_{n,k}$, where $y_2=y_1^2$, and $y_1$ is a first Stiefel Whitney class. Thus, $y_2$ is expressable using Stiefel Whitney classes over $V_{n,k}$. However, this does not allow us to deduce that $y_2$ is thus, expressible over $PV_{n,k}$. On the other hand if $n\equiv 0, 1 \pmod{4}$, the first possible value of $N$ is $4$, and if we further assume that $\binom{n}{4} \equiv 0 \pmod{2}$, then $y_3$ survives in the cohomology of $PV_{n,k}$. From \cite{Korbas2012}, we know that $y_3$ is not expressible in terms of Stiefel Whitney classes in $V_{n,k}$, so we have $\mbox{ucharrank}(PV_{n,k})=2$ in this case.  

\subsection{The complex case}
For the complex projective Stiefel manifold $PW_{n,k}=W_{n,k}/S^1$, recall that the $\Z_2$-cohomology is determined additively in \cite{Astey_Gitler1999}. 
$$ H^*(PW_{n,k};\Z_2)= \Z_2[x]/(x^N) \otimes \Lambda_{\Z_2}(A) , \text{ for } k<n ,$$ 
where  $N= \text{min}\{ j \mid n-k<j \leqslant n, \text{ and } \binom{n}{j} \text{ is odd } \} $, $ A= \{ y_j \mid n-k < j \leqslant n \} - \{ y_{_{N-1}}\} $ and $ |x|=2, |y_j|=2j-1 $. We compute the characteristic rank in the proposition below, which turns out to be much easier in this case. 
\begin{prop}
The upper characteristic rank of $PW_{n,k}$ is given by 
$$\mbox{ucharrank}(PW_{n,k})=\begin{cases} 2(n-k) & \mbox{if } 2\mid \binom{n}{k-1} \\ 
2(n-k)+2 &\mbox{if } 2\nmid \binom{n}{k-1}.\end{cases}$$ 
\end{prop}
\begin{proof}
Here the first class which appears in the exterior algebra part of the cohomology is in odd degree. The lower classes belong to the algebra generated by $x$ which is closed under products and Steenrod operations. By Wu's formula, a Stiefel Whitney class in degree which is not a power of $2$, is expressible in terms of lower degree Stiefel Whitney classes via the multiplication and Steenrod operations. This implies that the first generator of the exterior algebra part of the cohomology cannot be expressible in terms of Stiefel Whitney classes. The result follows. 
\end{proof}

\section{The circle quotient Stiefel manifolds}
There may be many circle actions defined on a real Stiefel manifold $V_{n,k}$. The ones that give us a homogeneous space act via a homomorphism  $S^1 \to O(k)$, where the latter acts on the Stiefel manifold by an orthogonal transformation on the vectors. This comes from a homomorphism into a maximal torus of $O(k)$, and up to conjugation they are a block diagonal inclusion of products of $SO(2)$. As in the case of projective Stiefel manifolds, we look at the diagonal inclusion of $S^1$ in the maximal torus. This gives us a manifold that we call the circle quotient Stiefel manifold, which is more precisely defined below.

We define the space $Y_{n,k}$ as the orbit space of the $S^1$-action on $V_{n,2k}$ defined as follows: We consider $S^1$ as $SO(2)$ embedded inside the maximal torus $SO(2)^k \subset SO(2k)$ by diagonal map. Then the action of our interest is the restriction of the action of $SO(2k)$ on $V_{n,2k}$ from right by matrix multiplication. \par
The construction gives a principal fiber bundle 
$$ S^1\longrightarrow V_{n,2k} \longrightarrow Y_{n,k}.$$ 
The complex line bundle associated to this principal bundle over $Y_{n,k}$ will be denoted by $\zeta$. We denote the realification of $\zeta$ by $\zeta_r$. \par 
From the definition of $Y_{n,k}$, it is clear that we have the following diagram of fibrations: 
\begin{equation}\label{univdiag}
    \xymatrix{
    	SO(2) \ar[r]\ar[d]_-\Delta & V_{n,2k} \ar[r] \ar@{=}[d] & Y_{n,k}\ar[d]_{\pi } \ar@/^1.5 pc/[dd]^q\\ 
SO(2)^k \ar[r]\ar[d] & V_{n,2k}\ar[r]\ar@{=}[d] & \widetilde{F}\ar[d]_p\\
O(2)^k\ar[r] & V_{n,2k}\ar[r] & F,
  }
\end{equation} where $\widetilde{F}$  is the $2^k$-sheeted cover of $F$, the space of flags $V_0\subset V_1\subset \cdots \subset V_k$ of subspaces of $\R^n$ with $\dim(V_i)=2i$. The canonical real vector bundles over $F$ will be denoted by $\xi_j$, for $1 \leqslant j \leqslant k+1$ and rank of $\xi_j=2$, for $1 \leqslant j \leqslant k$. Recall that the tangent bundle of $F$, $TF \cong \bigoplus_{1\leqslant i<j\leqslant k+1}\xi_i \otimes_\R\xi_j$, \cite{Lam1975}.
\subsection{The tangent bundle of $Y_{n,k}$}
 We now identify the tangent bundle of $Y_{n,k}$ in following theorem. 
\begin{thm}\label{tgtbdleY}
    The tangent bundle of $Y_{n,k}$ is described as  $$ TY_{n,k} = \binom{k}{2}\zeta_r\otimes_\R \zeta_r \oplus (k \zeta_r \otimes_\R q^*\xi_{k+1}) \oplus (k-1)\epsilon,$$ and additionally it satisfies the following relation $$TY_{n,k} \oplus \binom{k+1}{2}\zeta_r \otimes_\R \zeta_r = nk\zeta_r \oplus (k-1)\epsilon .$$ 
\end{thm}
\begin{proof}
From \ref{univdiag} we obtain the principal  bundle $$\frac{SO(2)^k}{\Delta(SO(2))} \longrightarrow Y_{n,k} 
\overset{\pi}\longrightarrow \widetilde{F}. $$ Now we can determine the tangent bundle of $Y_{n,k}$ as follows:
\begin{align*}
    TY_{n,k} &\cong \pi^*T\widetilde{F} \oplus (k-1)\epsilon\\
    &\cong q^* TF \oplus (k-1)\epsilon\\
    &\cong q^*\Big(\bigoplus_{1\leqslant i<j\leqslant k+1}\xi_i \otimes_\R\xi_j\Big)\oplus (k-1)\epsilon\\
    &\cong \binom{k}{2}\zeta_r\otimes_\R \zeta_r \oplus (k \zeta_r \otimes_\R q^*\xi_{k+1}) \oplus (k-1)\epsilon\\
 \end{align*}
Since $\bigoplus_{1 \leqslant j \leqslant k+1} \xi_j = n\epsilon$, we have $k\zeta_r \oplus q^*\xi_{k+1} =n\epsilon $. So we have $$TY_{n,k} \oplus \binom{k+1}{2}\zeta_r \otimes_\R \zeta_r = nk\zeta_r \oplus (k-1)\epsilon .$$
\end{proof}

Once we have the formula for the tangent bundle of $Y_{n,k}$, we may start analyzing questions regarding the number of linearly independent vector fields, restrictions on immersions into Euclidean space, skew embeddings, and characteristic rank. The necessary ingredient in the entire matter is the cohomology calculation for $Y_{n,k}$. 

\subsection{$\Z_2$-cohomology of $Y_{n,k}$}
As the action of $S^1$ on $V_{n,2k}$ occurs via a homomorphism $D : S^1 \to O(2k)$, we have the following commutative diagram 
\[ \xymatrix{S^1 \ar[r] \ar[d] & O(2k) \ar[d]\\
V_{n,2k} \ar[r] \ar[d] & V_{n,2k} \ar[d] \\
 Y_{n,k} \ar[r] \ar[d] & Gr_{2k}(\R^n)\ar[d] \\ 
                    \C P^\infty \ar[r] & BO(2k), }\]
where the rows are part of a homotopy fibration sequence. As a consequence the bottom square gives a homotopy pullback diagram 
\[ \xymatrix{ Y_{n,k} \ar[r] \ar[d] & Gr_{2k}(\R^n)\ar[d] \\ 
                    \C P^\infty \ar[r] & BO(2k). }\]
This implies that the circle quotient Stiefel manifold $Y_{n,2k}$ has the following universal property. 
\begin{prop}
Up to homotopy, the space $Y_{n,k}$ classifies complex line bundles $\xi$ such that for the realification $r(\xi)$, $kr(\xi)$ has a complimentary bundle $\mu$ of dimension $n-2k$. (That is, $\mu \oplus k r(\xi)=n \epsilon$.) 
\end{prop}

The fibrations above allow us to compute the cohomology of $Y_{n,k}$.
\begin{thm}\label{cohY}
	The $\Z_2$-cohomology of $Y_{n,k}$ is (additively)
	$$ H^*(Y_{n,k};\Z_2) \cong \Lambda_{\Z_2}(y_{n-2k}, \cdots,\widehat{y}_{_{2J-1}}, \cdots, y_{n-1}) \otimes \Z_2[x]/(x^J),  $$ where $deg(y_j)=j$ and $x$ is the mod $2$ Euler class of the bundle $\zeta$, and  $J=\text{min}\{ r \mid \binom{k+r-1}{k-1} \text{ is odd and } n-2k \leqslant 2r \leqslant n-1 \}$.
\end{thm}
\begin{proof}
    We consider the following commutative diagram of fibrations  
  \begin{equation}\label{cohYdiag}
  	\xymatrixcolsep{5pc}\xymatrix{V_{n,2k} \ar@{=}[d] \ar[r] & Y_{n,k} \ar[r] \ar[d] &\C P^\infty \ar[d] \\
  		V_{n,2k} \ar[r]  & Gr_{2k}(\R^n) \ar[r]  & BO(2k).
  	}
  \end{equation} 
It is known that $H^*(V_{n,k};\Z_2)$ is additively exterior algebra on generators $\{y_j \mid n-k \leqslant j <n \}$, with $deg(y_j)=j$, \cite{Borel1953}. We want to compute the action of differentials on $y_j$'s in the spectral sequence associated to first row of \ref{cohYdiag}. This is done by comparing the spectral sequences associated to the diagram \ref{cohYdiag}. But the spectral sequence for the bottom row of \ref{cohYdiag} was done in \cite{Basu2021} and says that $y_j$'s are transgressive with $\tau(y_j)= \bar{w}_{j+1}$, where $\bar{w}_{j+1}$ is the $(j+1)$-th inverse universal Stiefel Whitney class. So we know $y_j$'s are transgressive in the spectral sequence of our interest and  image of the transgression is given below.
    
    The map $\C P^\infty \longrightarrow BO(2k)$ is the classifying map for the $k$-fold Whitney sum of the underlying real $2$-plane bundle $r(\gamma^1)$ of the canonical complex line bundle $\gamma^1 $. So we can conclude $w(kr(\gamma^1))= (1+x)^k $ and hence $\overline{w}(kr(\gamma^1))= (1+x)^{-k} $. So \[ 
\tau(y_j)=\overline{w}_{j+1}(kr(\gamma^1)) = \begin{cases}  
 \binom{k+\frac{j+1}{2}-1}{k-1}x^{\frac{j+1}{2}}  & \mbox{if } j \mbox{ is odd} \\ 
 0 & \mbox{if } j \mbox{ is even}. 
 \end{cases} 
 \]
 This  calculation enables us to determine the $E_\infty$-page of the spectral sequence associated to the first row of the diagram and will be equal to $$E_\infty= V(S) \otimes \Z_2[x]/(x^J), $$ where $S= \{ y_j \mid n-2k \leqslant j \leqslant n-1 \}- \{ y_{_{2 J-1}}  \}$ and $J=\text{min}\{ r \mid \binom{k+r-1}{k-1} \text{ is odd and } n-2k \leqslant 2r \leqslant n-1 \}. $
 Now choosing lifts of the elements in $S$ we can conclude that additively $$ H^*(Y_{n,k};\Z_2) \cong \Lambda_{\Z_2}(y_{n-2k}, \cdots,\widehat{y}_{_{2J-1}}, \cdots, y_{n-1}) \otimes \Z_2[x]/(x^J),  $$ where $deg(y_j)=j$ and $x$ is the mod $2$ Euler class of the bundle $\zeta$, and  $J=\text{min}\{ r \mid \binom{k+r-1}{k-1} \text{ is odd and } n-2k \leqslant 2r \leqslant n-1 \}$.
 \end{proof}

\section{Computations for the circle quotient manifolds}
Computations described in section $3$ will allow us to obtain numerical information for certain topological invariants for spaces $Y_{n,k}$. 
\subsection{Stable span and paralleizability of $Y_{n,k}$}

We shall check the parallelizability of $Y_{n,k}$ using the relation obtained for $TY_{n,k}$ in \ref{tgtbdleY}.\par 
The total Pontryagin class for $\zeta_r$ is $p(\zeta_r)=1+x_0^2$, where $x_0 \in H^2(Y_{n,k};\Z) $ is the Euler class of $\zeta$ and whose mod $2$ reduction is $x$. To determine the total Pontryagin class for $\zeta_r \otimes_\R \zeta_r $, we note that the complexification of $\zeta_r \otimes_\R \zeta_r $ is $(\zeta\otimes_\C\zeta) \oplus (\zeta^*\otimes_\C\zeta^*) \oplus 2\epsilon $. Hence $p(\zeta_r \otimes_\R \zeta_r )=1+4x_0^2$. \par
So from the following relation relating the first Pontryagin classes: 
$$ p_1(TY_{n,k})= p_1(nk\zeta_r) - p_1(\binom{k+1}{2}\zeta_r \otimes_\R \zeta_r) ,$$
we get $p_1(TY_{n,k})=(nk-2k^2-2k)x_0^2$. Considering $n-2k \geq 4$, if we look at the Gysin sequence for the circle bundle $V_{n,2k} \longrightarrow Y_{n,k}$, we have 
$$ \xymatrix{0=H^3(V_{n,2k}) \ar[r] &H^2(Y_{n,k}) \ar[r]^{\cup x_0} &H^4(Y_{n,k})\ar[r] &\cdots }. $$ 
This ensures $x_0^2$ is non-zero and hence $p_1(TY_{n,k}) \neq 0$.
Hence we obtain 
\begin{thm}\label{Yparallel}
	If $n-2k \geqslant 4$, $Y_{n,k}$ is not parallelizable.
\end{thm} 
Recall that span of a vector bundle is its maximum number of linearly independent sections. We know Stiefel Whitney classes for a manifold $M$ provide an upper bound for its stable span, $\text{span}^0(M)=\text{span}(TM \oplus \epsilon) - 1$ . That bound for $Y_{n,k}$ is described in the theorem below. 
\begin{thm} \label{ISWY}
	If $$m := \text{max}\{ j \mid \binom{nk+j-1}{nk-1} \not\equiv 0  \textrm{ (mod}\ 2 ),  0 \leqslant j \leqslant J-1   \},$$ then $2m$-th inverse Stiefel Whitney class of $Y_{n,k}$, $\bar{w}_{2m}(Y_{n,k}) $ is non-zero and hence stable span of $Y_{n,k}$ satisfies the following inequality 
$$\text{span}^0(Y_{n,k}) \leqslant \text{dim}(Y_{n,k}) - 2m. $$
\end{thm}
\begin{proof}
	 The formula \ref{tgtbdleY} for tangent bundle of $Y_{n,k}$ allows us to calculate its total Stiefel Whitney class.
\begin{equation}\label{SWY}
	\begin{split}
	w(Y_{n,k}) &= w(nk\zeta_r)\cdot w\Big(\binom{k+1}{2}\zeta_r\otimes_{\R}\zeta_r\Big)^{-1}\\
               &= (1+x)^{nk} \cdot w\Big(\zeta_r\otimes_{\R}\zeta_r  \Big)^{-\binom{k+1}{2} }.
           \end{split}
\end{equation}
To determine the total Stiefel Whitney class of $\zeta_r\otimes_{\R}\zeta_r$, we apply splitting principle for this bundle. Suppose the bundle $\zeta_r$ splits as $L_1 \oplus L_2$ over $Y'$. Then $\zeta_r\otimes_{\R}\zeta_r$ splits as $\bigoplus_{1 \leqslant i,j \leqslant 2} L_i \otimes_{\R} L_j$ over $Y'$. We calculate the total Stiefel Whitney class of $\bigoplus_{1 \leqslant i,j \leqslant 2} L_i \otimes_{\R} L_j$ below:
\begin{align*}
	w\big(\bigoplus_{1 \leqslant i,j \leqslant 2} L_i \otimes_{\R} L_j \big) &= \prod_{1 \leqslant i,j \leqslant 2} w(L_i \otimes_{\R} L_j)\\ 
	&= \prod_{1 \leqslant i,j \leqslant 2} (1+w_1(L_i)+w_1(L_j))\\
	&= (1+w_1(L_1)+w_1(L_2))^2\\
	&=(1+w_1(L_1 \oplus L_2))^2=1.
\end{align*} 
Hence we must have  $w\Big(\zeta_r\otimes_{\R}\zeta_r  \Big)=1$. So from \ref{SWY}, we get
\begin{equation}\label{SWY1} 
	w(Y_{n,k})=(1+x)^{nk} \implies \bar{w}(Y_{n,k})=(1+x)^{-nk} .
\end{equation} 
Then from the definition of $m$, we see that $\bar{w}_{_{2m}}(Y_{n,k}) \neq 0$ and the theorem follows.
\end{proof}

\subsection{Skew embedding and immersion dimensions of $Y_{n,k}$} 
The result concerning non-vanishing of inverse Stiefel Whitney class stated in theorem \ref{ISWY} immediately produces lower bounds of skew embedding and immersion dimension of $Y_{n,k}$. 
\begin{prop}
    $Y_{n,k}$ does not admit an immersion in $\R^{dim(Y_{n,k})+2m-1}$ and a skew embedding in $\R^{2 dim(Y_{n,k}) + 4m}$, where $$m = \text{max}\{ j \mid \binom{nk+j-1}{nk-1} \not\equiv 0  \textrm{ (mod}\ 2 ), 0 \leq j \leq J-1   \}.$$
\end{prop}
\begin{proof}
    The statement concerning immersion follows since Stiefel Whitney classes of the normal bundle of an immersion in an Euclidean space are same as the inverse Stiefel Whitney classes of tangent bundle.\par 
    The statement about skew embedding follows directly by combining theorems \ref{skewCond} and \ref{ISWY}.
\end{proof}
\subsection{Characteristic rank of $Y_{n,k}$} We shall determine the upper characteristic rank of $Y_{n,k}$ for a large number of cases.
\begin{thm}
If $n-2k=5$, $6$ or $\geq 9$, the upper characteristic rank of $Y_{n,k}$ is given by 
$$\mbox{ucharrank}(Y_{n,k})=\begin{cases} n-2k-1 & \mbox{if } 2\mid \binom{n}{2k-1} \\ 
	n-2k &\mbox{if } 2\nmid \binom{n}{2k-1}.\end{cases}$$    
\end{thm}
\begin{proof} The proof is similar to the proof of theorem \ref{charrkPV}.
    First we consider the case  $n-2k \neq2J-1$. Then there is a non-zero class $y_{n-2k}\in H^{n-2k}(Y_{n,k})$ which is pulled back to the generator of $H^{n-2k}(V_{n,2k};\Z_2)=\Z_2$ and which is not expressible as a polynomial in Stiefel Whitney classes of some bundle over $Y_{n,k}$ because that would contradict the fact that $\text{ucharrank}(V_{n,2k})= n-2k-1 $, \cite{Korbas2012} otherwise. So in this case we must have $\text{ucharrank}(Y_{n,k})=n-2k-1 $ since $ \text{charrank}(\zeta)=n-2k-1 $.\par  
    Next we consider the case $n-2k=2J-1$. Then $\bigoplus_{0\leqslant i \leqslant n-2k} H^i(Y_{n,k};\Z_2)= \Z_2[x]/(x^J)$ and $H^{n-2k+1}(Y_{n,k};\Z_2)=\Z_2\{y_{n-2k+1}\}$. So for $n-2k+1\neq$ a power of $2$, Wu's formula rules out the possibility of $y_{n-2k+1}$ being Stiefel Whitney class of a bundle over $Y_{n,k}$ and we get $\text{ucharrank}(Y_{n,k})=n-2k$. And for $n-2k+1 =$ a power of 2 greater than $8$, invoking the lemma \ref{lempvnk} one again obtains $\text{ucharrank}(Y_{n,k})=n-2k$.
    
\end{proof}

\bibliographystyle{siam}
\bibliography{my_bibliography}
\clearpage

\end{document}